\documentclass[preprint,11pt]{imsart}
\RequirePackage[T1]{fontenc}
\usepackage{color}
\usepackage[latin1]{inputenc}
\usepackage[T1]{fontenc}

\usepackage{amsmath,amssymb,graphicx}
\usepackage{epstopdf}
\usepackage{color}
\RequirePackage{amsfonts,amsthm,amsmath}
\RequirePackage[colorlinks,citecolor=blue,urlcolor=blue]{hyperref}
\usepackage{a4wide}
\usepackage{here}
\usepackage{float}
\usepackage{amsmath}
\usepackage{bbm}
\allowdisplaybreaks

\newcommand{\nbR}{\mathbb{R}}
\newcommand{\nbN}{\mathbb{N}}
\newcommand{\nbQ}{\mathbb{Q}}

\newcommand{\nbu}{\mathbbm{1}}

\newcommand{\nbP}{\mathbb{P}}

\newcommand{\nbE}{\mathbb{E}}

\newtheorem{theorem}{{\bf Theorem}}
\newtheorem{lemm}[theorem]{{\bf Lemma}}
\newtheorem{cor}[theorem]{{\bf Corollary}}

\newtheorem{prop}[theorem]{\bf Proposition}

\def\PP{{\mathbb P}}

\newcommand{\N}{\ensuremath{\mathbb{N}}}

\begin{document}
\linespread{1.5}
\begin{frontmatter}

\title{Asymptotics for infinite server queues with fast/slow Markov switching and fat tailed service times}

\runtitle{Infinite server queues with switching and fat tailed service times}

\begin{aug}
\author{\fnms{\Large{Landy}}
\snm{\Large{Rabehasaina}}
\ead[label=e2]{lrabehas@univ-fcomte.fr}} 
\runauthor{L.Rabehasaina}

\address{\hspace*{0cm}\\
Laboratory of Mathematics, University Bourgogne Franche-Comt\'e,\\
16 route de Gray, 25030 Besan\c con cedex, France.\\[0.2cm]
\printead{e2}}

\end{aug}

\vspace{0.5cm}

\begin{abstract}
We study a general $k$ dimensional infinite server queues process with Markov switching, Poisson arrivals and where the service times are fat tailed with index $\alpha\in (0,1)$. When the arrival rate is sped up by a factor $n^\gamma$, the transition probabilities of the underlying Markov chain are divided by $n^\gamma$ and the service times are divided by $n$, we identify two regimes ("fast arrivals", when $\gamma>\alpha$, and "equilibrium", when $\gamma=\alpha$) in which we prove that a properly rescaled process converges pointwise in distribution to some limiting process. In a third "slow arrivals" regime, $\gamma<\alpha$, we show the convergence of the two first joint moments of the rescaled process.
\end{abstract}

\end{frontmatter}

\normalsize


\section{Introduction and notation}\label{sec:model}
\subsection{Model and related work}
The classical infinite server queue consists of a system where tasks or customers arrive according to a general arrival process and begin receiving service immediately. Such a model was studied extensively, under various assumptions on the interarrival and service time distributions, in \cite[Chapter 3, Section 3]{T62}. Several variants or extensions have been considered, in particular where arrivals and service times are governed by an external background Markovian process \cite{OP86, D08, FA09, BKMT14, MDT16, BDTM17}, or where customers arrive in batches \cite{MT02}. An extension to a network of infinite-server queues where arrival and service rates are Markov modulated is studied in \cite{JMDW19}. 

We consider yet another generalization of this model with Markov switching described as follows. Let $\{ N_t,\ t\ge 0\}$ be a Poisson process with intensity $\lambda>0$, corresponding jump times $(T_i)_{i\in\nbN}$ satisfying $T_0=0$, such that $(T_i-T_{i-1})_{i\ge 1}$ is a sequence of independent and identically distributed (iid) random variables with same exponential distribution with parameter $\lambda>0$, denoted by ${\cal E(\lambda)}$. Let $(L_{ij})_{i\in \nbN, j=1,\ldots,k}$ be a sequence of independent random variables such that the sequence of vectors $\left((L_{i1},\ldots,L_{ik})\right)_{i\in \nbN}$ is iid (with entries $L_{i1}$,\ldots,$L_{ik}$ having different distributions for each $i$). Finally, for some $K$ and $k$ in $\nbN^*$ we consider the discrete set ${\cal S}=\{ 0,\ldots,K\}^k$ and a stationary finite Markov chain $(X_i)_{i\in\nbN}$ with state space ${\cal S}$. Then, for all $i$, $X_i$ is a vector of the form $X_i=(X_{i1},\ldots,X_{ik})$ with $X_{ij}\in \{ 0,\ldots,K\}$, $j=1,\ldots,k$. We then define the following $k$ dimensional process $\{Z(t)=(Z_1(t),\ldots,Z_k(t)),\ t\ge 0 \}$ with values in $\nbN^k$ as
\begin{equation}\label{def_Z_t}
Z_j(t)=\sum_{i=1}^{N_t}X_{ij} \nbu_{[t<L_{ij}+T_i]}=\sum_{i=1}^{\infty}X_{ij} \nbu_{[T_i\le t<L_{ij}+T_i]},\quad j=1,...,k.
\end{equation}
The process defined by \eqref{def_Z_t} has many applications, of which we list two most important ones:
\begin{itemize}
\item {\it incurred but not reported correlated claims: } in an actuarial context, $Z(t)=(Z_1(t),\ldots,Z_k(t))$ represents a set of branches where $Z_j(t)$ is the number of incurred but non reported (IBNR) claims in the $j$th branch of an insurance company. Here $X_{ij}$ is the number of such claims arriving in that branch at time $T_i$, and $L_{ij}$ is the related delay time before the claim $j$ is reported. From another point of view, $X_{ij}\in [0,+\infty)$ may also represent the amount (say, in euros) of the claim occurring at time $T_i$ in the $j$th branch, in which case $Z_j(t)$ is the total amount of undeclared claims which have occurred by time $t$. Another application is when $K=1$, in which case $X_{ij}=0$ means that the claim in branch $j$ occurring at time $T_i$ is reported and dealt with immediately by the policyholder, whereas $X_{ij}=1$ means that some effective lag in the report is observed. The Markovian nature of $(X_i)_{i\in\N}$ here is important from a practical point of view, as a claim amount at time $T_i$ may impact the one at time $T_{i+1}$, or because a policyholder may decide to grant a long report delay for the claim at time $T_{i+1}$ with high probability if the claim at time $T_i$ is reported immediately.
\item {\it infinite server queues with batch arrivals and Markov switching: } $Z(t)=(Z_1(t),\ldots,Z_k(t))$ represents a set of $k$ correlated queues with an infinite number of servers, such that customers arrive at each time $T_i$, with $X_{ij}$ customers arriving in queue $j\in\{1,...,k\}$, with corresponding (same) service times $L_{ij}$ (as an example, the basic case where $X_{ij}=1$ for all $i\in \nbN$ and $j=1,...,k$ corresponds to $k$ customers arriving simultaneously at each instant $T_i$). $Z_j(t)$ can also be seen as the number of customers of class $j$ in a (single) infinite-server queue, as illustrated in \cite[Figure 1]{RW16}. Other infinite-server queue, such as one where the customers within a batch arriving at time $T_i$ have different service times, may be inferred from the model \eqref{def_Z_t} by choosing an appropriate value of $k$ and Markov chain $(X_i)_{i\in\N}$, see \cite[Section 6]{RW16}. Here, the Markov switching is a major novelty in the present model because it allows for some dependence between the successive number of incoming customers. One simple example is when $K=1$, so that $X_{ij}=0$ means that an incoming customer at time $T_i$ in queue $j$ is rejected from the system, whereas $X_{ij}=1$ means that it is accepted: a classical situation would then be that if a customer is rejected at time $T_i$ then the next one could be accepted with high probability, at time $T_{i+1}$. In other words, this Markov switching can help model traffic regulation mechanisms.
\end{itemize}
The present paper follows \cite{RW18}, which studies the transient or limiting distribution of a discounted version of $Z(t)$ of the form
\begin{equation}\label{discounted_Z(t)}
Z_j(t)=\sum_{i=1}^{N_t}X_{ij} e^{-a (L_{ij}+T_i)}\nbu_{[t<L_{ij}+T_i]}=\sum_{i=1}^{\infty}X_{ij} e^{-a (L_{ij}+T_i)} \nbu_{[T_i\le t<L_{ij}+T_i]},\quad j=1,...,k,
\end{equation}
for $t\ge 0$. The main difference with \cite{RW18} is that the latter has more general assumptions on the interarrival and service distributions, whereas we focus here on Poisson arrivals. Even though the assumptions are more restrictive than in \cite{RW18}, the goal here is different in that we are trying to exhibit different behaviours for the limiting models when the arrival rate is increased and the service times are decreased by suitable factors, whereas \cite{RW18} is more focused on analytical stochastic properties such as the moments of $Z(t)$ or its limiting distribution as $t\to\infty$. The discounting factor $a \ge 0$ in \eqref{discounted_Z(t)} is important in situations e.g. where, in an actuarial context, $X_{ij} e^{-a (L_{ij}+T_i)}$ represents the value of the claim amount at the actual realization time $L_{ij}+T_i$.
Furthermore, the state space ${\cal S}=\{ 0,\ldots,K\}^k$, although seemingly artificially complex, allows in fact for some flexibility and enables us to retrieve some known models. In particular, consider a Markov-modulated infinite-server queue, i.e. a queueing process $\{{\cal Z}(t),\ t\ge 0\}$ of which interarrivals and service times are modulated by a background continuous time Markov chain $\{Y(t),\ t\ge 0\}$ with state space say $\{1,...,\kappa\}$, i.e. such that customers arrive on the switching times of the Markov chain, with service times depending on the state of the background process (see \cite{MT02}, \cite[Model II]{MDT16}). Then \cite[Section 6]{RW18} explains how this process $\{{\cal Z}(t),\ t\ge 0\}$ can be embedded in a process $\{Z(t),\ t\ge 0\}$ defined by \eqref{def_Z_t} with an appropriate choice of $k$, $K$ in function of $\kappa$, as well as of the Markov chain $(X_i)_{i\in\nbN}$ and the sequence $(L_{ij})_{i\in \nbN, j=1,\ldots,k}$ of service times. Thus, studying a general process $\{Z(t),\ t\ge 0\}$ in \eqref{def_Z_t} allows to study a broad class of infinite server queue models in a similar Markov modulated context. 

We now proceed with some notation related to the model and used throughout the paper. Let $P=(p(x,x'))_{(x,x')\in {\cal S}^2}$ and $\pi=(\pi(x))_{x\in {\cal S}}$ (written as a row vector) be respectively the transition matrix and stationary distribution of the Markov chain. We next define for all $r\ge 0$ and $s=(s_1,\ldots,s_k)\in (-\infty,0]^k $,
\begin{align}
\tilde{\pi}(s,r)&:= \displaystyle\mbox{diag}\left[ \nbE \left( \exp\left\{ \sum_{j=1}^k s_jx_j  \nbu_{[L_j>r]}\right\}\right),\ x=(x_1,\ldots,x_k)\in{\cal S}\right],\label{def_pi_Q_tilda}\\
\Delta_i&:=\mbox{diag} \left[ x_i,\ x=(x_1,\ldots,x_k)\in{\cal S}\right],\quad i=1,\ldots,k,\label{Di}
\end{align}
where $P'$ denotes the transpose of matrix $P$. $I$ is the identity matrix, ${\bf 0}$ is a column vector with zeroes, and ${\bf 1}$ is a column vector with $1$'s, of appropriate dimensions.
The Laplace Transform (LT) of the process $Z(t)$ jointly to the state of $X_{N_t}$ given  the initial state of $X_0$ is denoted by
\begin{equation}\label{def_mgf}
\psi(s,t):=\left[ \nbE\left( \left. e^{<s,Z(t)>}\nbu_{[X_{N_t}=y]}\right| X_0=x\right)\right]_{(x,y)\in {\cal S}^2},\quad t\ge 0,\ s=(s_1,\ldots,s_k)\in (-\infty,0]^k
\end{equation}
where $<\cdot, \cdot>$ denotes the Euclidian inner product on $\nbR^k$. Note that $X_0$ has no direct physical interpretation here, as the claims sizes/customer batches are given by $X_i$, $i\ge 1$, and is rather introduced for technical purpose.

We finish this section with the following notation. For two sequences of random variables $(A_n)_{n\in\nbN}$ and $(B_n)_{n\in\nbN}$ and two random variables $A$ and $B$, the notation ${\cal D}\left(\left. A_n\right|B_n\right)\longrightarrow_{n\to\infty} {\cal D}\left(\left. A\right|B\right)$ indicates that, as $n\to \infty$, the conditional distribution of $A_n$ given $B_n$ converges weakly to the conditional distribution of $A$ given $B$.


\subsection{Rescaling}\label{sec:rescale}
We arrive at the main topic of the paper, which is to be able to provide some information on the distribution of $Z(t)$ in \eqref{def_Z_t}. In the particular case of Poisson arrivals, and since $Z(t)$ in \eqref{def_Z_t} is a particular case of the process in \eqref{discounted_Z(t)} with discount factor $a=0$, the LT $\psi(s,t)$ defined in \eqref{def_mgf} is characterized by  \cite[Proposition 4]{RW18}, which we rewrite here:
\begin{prop}\label{prop_Poisson_psi}
When $\{ N_t,\ t\ge 0\}$ is a Poisson process with intensity $\lambda>0$, then $\psi(s,t)$ is the unique solution to the first order linear (matrix) differential equation
\begin{equation}\label{Poisson_ODE}
\partial_t \psi(s,t) =[\lambda (P-I) + \lambda P(\tilde{\pi}(s,t)-I)]\psi(s,t)
\end{equation}
with the initial condition $\psi(s,0)= I$.
\end{prop}
Unfortunately, the first order ordinary differential equation \eqref{Poisson_ODE} does not have an explicit expression in general, so that studying the (transient or stationary) distribution of the couple $(Z(t),X_{N_t})$ is difficult. In that case, as in \cite{BDTM17, BKMT14, MDT16}, it is appealing to study the process when the intensity of the Poisson process is sped up and the switching rates of the Markov chain are modified. Similarly to those papers, the goal of this paper is thus to study the behaviour of the queue/IBNR process in "extreme conditions" for the arrival rates, transition rates and delays, while trying to maintain minimal assumptions on the service time distributions. For this we will suppose that the rescalings of the parameters, denoted by ${\bf (S1)}$, ${\bf (S2)}$ and ${\bf (S3)}$ hereafter, are performed as follows:
\begin{itemize}
\item the arrival rate is multiplied by $n^\gamma$ for some $\gamma>0$, denoted by $${\bf (S1)}\quad \lambda_n= \lambda n^\gamma ,$$
with associated Poisson process $\{ N_t^{(n)},\ t\ge 0\}$ and jump times $(T_i^n)_{i\in\nbN}$,
\item the transition probabilities $p(x,y)$ are divided by $n^\gamma$ when $x\neq y$, $x$, $y$ in ${\cal S}$, i.e. the new transition matrix is given by $${\bf (S2)} \quad P_n=P/n^\gamma + (1-1/n^\gamma) I,$$
with corresponding stationary Markov chain $(X_i^{(n)})_{i\in\nbN}$, having the same distribution $\pi$ as $(X_i)_{i\in\nbN}$.
\end{itemize}
Since the transition matrix $P_n$ verifies $P_n\longrightarrow_{n\to\infty} I$, such normalizing assumptions imply that, as $n\to\infty$,  one is close to a model where the arriving customers or claims come in the $k$ queues in batches with same fixed size: those queues are nonetheless correlated because the customers arrive according to the same Poisson process. Also, observe that
$\lambda_n (P_n-I)$ is the infinitesimal generator of the continuous time Markov chain $\left\{Y^{(n)}(t)=X^{(n)}_{N^{(n)}_t},\ t\ge 0\right\}$ of which embedded Markov chain is the underlying Markov chain i.e. $(Y^{(n)}(T_i^n))_{i\in\nbN}=(X_i^{(n)})_{i\in\nbN}$. Thus, since the rescaling is such that
$$\lambda_n (P_n-I)=\lambda (P-I)$$
for all $n$ (a property which will be extensively used in the paper), we remark that the rescalings {\bf (S1)} and {\bf (S2)} for the
arrival rate and the transition probabilities are such that the transition rates between the states of $\cal S$ of $\{Y^{(n)}(t),\ t\ge 0\}$ are independent from $n$, which allows for enough dynamics in the model that compensates the fact that $P_n$ tends to $I$, and yielding non trivial asymptotics in the convergence results in this paper as $n\to\infty$.

The assumptions for the service times/delays distribution are the following. We first suppose that the base model features fat tailed distributed service times with same index $\alpha \in (0,1)$, i.e. such that
$$\nbP(L_j >t)\sim 1/t^\alpha,\quad t\to \infty,$$
for all $j=1,...,k$. This kind of distribution (included in the wider class of heavy tailed distributions) mean that the service times are "large". In particular, those service times have {\it infinite expectation}. Furthermore, the rescaling for the service times is such that they are divided by $n$, denoted by
$$
{\bf (S3)}\quad L_j^{(n)}=L_j/n.
$$ 
Hence, the situation is the following: the arrivals are sped up by factor $n^\gamma$, but this is compensated by the fact that the delay times are diminished with factor $n$, so that one expects one of the three phenomena to occur at time $t$ for the limiting model: the arrivals occur faster than it takes time for customers to be served and the corresponding queue content $Z^{(n)}(t)$ grows large as $n\to\infty$, the arrivals occur slower and services are completed fast so that $Z^{(n)}(t)$ tends to $0$ as $n\to\infty$, or an equilibrium is reached. Those three cases will be studied in the forthcoming sections. Some limiting behaviour was studied in \cite{BKMT14, MDT16}, where the authors identified three regimes for different scalings in a Markov modulating context and obtained a Central Limit Theorem for a renormalized process, 
when the service times have general distribution with finite expectation or are exponentially distributed. \cite{BDTM17} provides some precise asymptotics on the tail probability of the queue content for exponentially distributed service times. \cite{JMDW19} provides a diffusion approximation for a model with exponentially distributed service times. A novelty in this paper is that we restrict here the class of distributions to that of fat tailed distributions in order to exhibit (under different scalings) a different behaviour and different limiting distribution which is not gaussian. Also note that the class of fat tailed distributions is interesting in itself as, in actuarial practice, this corresponds to {\it latent claims}, i.e. very long delays which are incidentally in practice often not observed (as the case $\alpha\in(0,1)$ corresponds to the $L_j$'s having infinite expectation), see \cite[Section 6.6.1]{H17}. This motivates the convergence results in this paper, which feature the exponent $\alpha$ as the only information required on those delays. This in itself is a noticeable difference from the Central Limit Theorems obtained in \cite[Section 4]{BKMT14}, where the normalization and limiting distribution require the explicit cumulative distribution function of the service times. Not only that, but the scaling is rather done in those references \cite{BKMT14, MDT16} on the transition rates of the underlying continuous time Markov chain modulating the arrival and service rates, whereas here these are constant, as we saw that $\lambda_n (P_n-I)=\lambda (P-I)$ is independent of $n$, and the scaling is rather done on the service times instead. When the service times are heavy tailed, this particular model can also be seen as a generalization of the {\it infinite source model}, see \cite[Section 2.2]{MRRS02}. Since the class of fat tailed distributions is a sub-class of the set of heavy tailed distributions, the normalizations {\bf (S1)} and {\bf (S3)} can be directly compared to \cite[Section 3.1]{MRRS02}, which studies limiting distributions of such normalized processes, and where the authors introduce the notion of so-called Slow and Fast Growth conditions when the arrival rate of customers is respectively negligible or dominant, compared to the service times. The reader is also referred to \cite{GK03} for a a similar model where the interarrivals are heavy tailed.
All in all, what is going to be studied hereafter is, when $t$ is fixed in say $[0,1]$ w.l.o.g., the limiting distribution as $n\to\infty$ of the $\nbN^k\times {\cal S}$ valued random vector
$$
\left(Z^{(n)}(t),X^{(n)}_{N^{(n)}_t}\right)
$$
under rescaling ${\bf (S1)}$, ${\bf (S2)}$ and ${\bf (S3)}$, or of a renormalized version of it in the "fast" or "slow" arriving customers case. Note that that the convergence is proved on the interval $[0,1]$, but all proofs can be adapted to show the convergence on any interval $[0,M]$ for $M>0$. The corresponding joint Laplace Transform is given by
\begin{equation}\label{def_LT_n}
\psi^{(n)}(s,t)=\left[ \nbE\left( \left. e^{<s,Z^{(n)}(t)>}\nbu_{\left[X^{(n)}_{N^{(n)}_t}=y\right]}\right| \ X_0^{(n)}=x\right)\right]_{(x,y)\in {\cal S}^2},\ s=(s_1,...,s_j)\in (-\infty, 0]^k,
\end{equation}
where we recall that $(X^{(n)}_i)_{i\in \nbN}$ is the underlying Markov chain with generating matrix $P_n$, stationary distribution $\pi$, and $\left\{N^{(n)}_t,\ t\ge 0\right\}$ is a Poisson process representing the arrivals, with scaled intensity $\lambda_n$. We also introduce the first and second joint matrix moments defined by
\begin{equation}\label{def_moments}
\begin{array}{rcl}
M_j^{(n)}(t)&:= &\left[ \nbE\left(\left. Z^{(n)}_j(t) \nbu_{\left[X^{(n)}_{N^{(n)}_t}=y\right]}\right| \ X_0^{(n)}=x\right) \right]_{(x,y)\in {\cal S}^2},\quad j=1,...,k,\\
M_{jj'}^{(n)}(t)&:= &\left[ \nbE\left(\left. Z^{(n)}_j(t)\ Z^{(n)}_{j'}(t) \nbu_{\left[X^{(n)}_{N^{(n)}_t}=y\right]}\right| \ X_0^{(n)}=x\right) \right]_{(x,y)\in {\cal S}^2},\quad j,j'=1,...,k .
\end{array}
\end{equation}
\section{Statement of results and organization of paper}
The core results of the paper concerning the different regimes are given in the following two Theorems \ref{theo_regimes} and \ref{theo_slow_arrival}:
\begin{theorem}\label{theo_regimes}
Let $\{{\cal X}^\alpha (t)=({\cal X}^\alpha_1 (t),...,{\cal X}^\alpha_k (t)),\ t\in[0,1] \}$ be a $\{0,...,K\}^k$ valued continuous time inhomogeneous Markov chain with infinitesimal generating matrix $\frac{1}{1-\alpha} (1-t)^{\frac{\alpha}{1-\alpha}}\lambda (P-I)$ with ${\cal X}^\alpha (0)\sim \pi$, and $\{\nu_j^\alpha(t),\ t\in [0,1] \}$, $j=1,...,k$, be $k$ independent Poisson processes with same intensity $\frac{\lambda}{1-\alpha}$, independent from $\{{\cal X}^\alpha (t),\ t\in[0,1] \}$. Let $t\in[0,1]$ fixed.
\begin{itemize}
\item {\bf Fast arrivals: }If $ \gamma >\alpha$ then, as $n\to\infty$,
\begin{multline}\label{convergence_fast}
{\cal D}\left(\left. \left(\frac{Z^{(n)}(t)}{n^{\gamma-\alpha}},X^{(n)}_{N^{(n)}_t} \right)\right|\ X^{(n)}_0\right)\\
\longrightarrow {\cal D}\left(
\left. 
\left(\frac{\lambda}{1-\alpha} \int_{1-t^{1-\alpha}}^{1} {\cal X}^\alpha(v)\ dv,\ {\cal X}^\alpha (1)\right)
\right|\ {\cal X}^\alpha\left(1-t^{1-\alpha}\right)
\right),
\end{multline}
\item {\bf Equilibrium: } If $\gamma =\alpha$ then, as $n\to\infty$,
\begin{multline}\label{convergence_equilibrium}
{\cal D}\left(\left.  \left(Z^{(n)}(t),X^{(n)}_{N^{(n)}_t} \right)\right|\ X^{(n)}_0\right)\\
\longrightarrow {\cal D}\left(
\left. 
\left( \left(\int_{1-t^{1-\alpha}}^{1} {\cal X}^\alpha_j(v)\ \nu_j^\alpha(dv)\right)_{j=1,...,k},\ {\cal X}^\alpha (1)\right)
\right|\ {\cal X}^\alpha\left(1-t^{1-\alpha}\right)
\right).
\end{multline}
\end{itemize}
\end{theorem}
We note that the terms in the limits on the right hand side of \eqref{convergence_fast} and \eqref{convergence_equilibrium} feature simple objects (in regards to the complexity of the original model) where the only characteristic parameters needed are $\lambda$, $P$ and $\alpha$; in particular, and apart from $\alpha$, characteristics of the service times $L_j$, $j=1,...,k$, such as their cumulative distribution functions, do not show up in the limiting distributions \eqref{convergence_fast} and \eqref{convergence_equilibrium}.
The convergences in distribution \eqref{convergence_fast} and \eqref{convergence_equilibrium} give some information on convergence of the (possibly renormalized) joint distribution of the couple $\left(Z^{(n)}(t),X^{(n)}_{N^{(n)}_t}\right)$, $t\in [0,1]$ given the initial state $X_0^{(n)}$. 
Intuitively, for fixed $t\in [0,1]$, in the right hand sides of \eqref{convergence_fast} and \eqref{convergence_equilibrium} we may interpret the inhomogeneous continuous time Markov chain $\{{\cal X}^\alpha (v),\ v\in[1-t^{1-\alpha},1] \}$ as the limiting counterpart of the modulating process $\left\{ X^{(n)}_{N^{(n)}_v},\ v\in [0,t]\right\}$. On an even cruder level, we observe in the Fast arrivals case from \eqref{convergence_fast} that each entry of $Z^{(n)}(t)=(Z^{(n)}_1(t),...,Z^{(n)}_k(t))$ behaves roughly like $n^{\gamma-\alpha}t^{1-\alpha}$. The intuition behind this behaviour may be explained as follows. Within queue $j=1,...k$, there are approximately $\lambda n^{\gamma} t$ arrivals in the interval $[0,t]$, each arriving customer with service time distributed as $L^{(n)}_j$, so that we may very grossly consider that a customer is still present at time $t$ with probability $\nbP (L^{(n)}_j>t)=\nbP (L_j/n>t) $. Hence the number of customers in queue $j$ is approximately
$$
Z^{(n)}_j(t)\approx \lambda n^{\gamma} t \times \nbP (L^{(j)}/n>t) = \lambda n^{\gamma} t \times\nbP (L^{(j)}>nt) \approx \lambda n^{\gamma} t \times \frac{1}{(nt)^\alpha}=\lambda n^{\gamma-\alpha}t^{1-\alpha}
$$
which is the expected order of growth $n^{\gamma-\alpha}t^{1-\alpha}$. Of course, such approximations are very crude, however this enables us to justify the presence of the normalizing factor $n^{\gamma-\alpha}$ as well as the time dilated factor $t^{1-\alpha}$ in \eqref{convergence_fast}.

In the case when $ \gamma <\alpha$, proving the convergence in distribution of an adequate normalization of $Z^{(n)}(t)$ seems more difficult. The following result show that the two first moments of $Z^{(n)}(t)$ converge under respective normalization $n^{\alpha-\gamma}$ and $n^{(\alpha-\gamma)/2}$:
\begin{theorem}[\bf Slow arrivals]\label{theo_slow_arrival}
If $ \gamma <\alpha$ then the following convergences of the two joint moments hold as $n\to \infty$
\begin{eqnarray}
n^{\alpha-\gamma} M_j^{(n)}(t)&\longrightarrow & \lambda e^{\lambda t (P-I)}\int_0^t \frac{1}{v^\alpha} e^{-\lambda v (P-I)} \Delta_j e^{\lambda v (P-I)} dv, \label{convergence_slow_M1}\\
n^{\alpha-\gamma} M_{jj}^{(n)}(t)&\longrightarrow &  \lambda e^{\lambda t (P-I)}\int_0^t \frac{1}{v^\alpha} e^{-\lambda v (P-I)} \Delta_j^2 e^{\lambda v (P-I)} dv, \label{convergence_slow_M20}\\
n^{\alpha-\gamma} M_{jj'}^{(n)}(t)&\longrightarrow &  0 \quad j\neq j', \label{convergence_slow_M21}
\end{eqnarray}
for all $j$, $j'\neq j$, in $1,...,k$, $t\in [0,1]$, where we recall that $\Delta_j$ is defined in \eqref{Di}.
\end{theorem}
One interesting by-product of Theorem \ref{theo_regimes} is that it in particular gives some insight on the (non conditional) limiting distribution of $Z^{(n)}(t)$, with a limit in a somewhat simpler form. More precisely, the following corollary follows from the proofs of \eqref{convergence_fast} and \eqref{convergence_equilibrium}:
\begin{cor}\label{rem:marginal}
In the Fast arrivals case $\gamma >\alpha$, the following convergence holds
\begin{equation}\label{remark_conv_distrib_simpler_fast}
\frac{Z^{(n)}(t)}{n^{\gamma-\alpha}}\stackrel{\cal D}{\longrightarrow} \lambda \int_{0}^{t} \frac{{\cal Y}(v)}{v^\alpha}\ dv,\ n\to\infty, \quad t\in [0,1], 
\end{equation}
where $\{{\cal Y}(t)=({\cal Y}_1(t),...,{\cal Y}_k(t)),\ t\in [0,1]\}$ is a (time homogeneous) stationary continuous time Markov chain on the state space ${\cal S}$, with infinitesimal generator matrix defined by
\begin{equation}\label{generator_Y}
 \lambda (\Delta_\pi^{-1} P' \Delta_\pi-I),\quad \Delta_\pi:=\mbox{\normalfont diag}(\pi(x),\ x\in {\cal S}).
\end{equation}
In the Equilibrium case $\gamma =\alpha$, one has
\begin{equation}\label{remark_conv_distrib_simpler_equilibirum}
Z^{(n)}(t)\stackrel{\cal D}{\longrightarrow}\left(  \int_0^t {\cal Y}_j(v) \ \tilde{\nu}_j^\alpha(dv)\right)_{j=1,...,k},\ n\to\infty, \quad t\in [0,1].
\end{equation}
Here $\{\tilde{\nu}_j^\alpha(t),\ t\in [0,1] \}$, $j=1,...,k$, are the inhomogeneous Poisson processes defined by $\tilde{\nu}_j^\alpha(t)=\nu_j^\alpha(t^{1-\alpha})$, $t\in [0,1]$, where $\{\nu_j^\alpha(t),\ t\in [0,1] \}$, $j=1,...,k$, are defined in Theorem \ref{theo_regimes}.
\end{cor}
As mentioned in Section \ref{sec:rescale}, \cite{MRRS02, GK03} introduced a notion of Fast and Slow growth similar to the Fast and Slow arrivals presented in Theorems \ref{theo_regimes} and \ref{theo_slow_arrival}, for a process of interest which is either a superposition of renewal processes with heavy tailed interarrivals or the cumulative input of an infinite source Poisson model with heavy tailed services. In those references, the process is shown to converge weakly or in finite dimensional distributions towards specific limit processes under appropriate scaling, see \cite[Theorem 1]{MRRS02} and \cite[Theorem 1]{GK03}. Here, the outline of the proof of Theorem \ref{theo_regimes} is the following:
\begin{itemize}
\item We will first expand the LT of the left hand side of \eqref{convergence_fast} and \eqref{convergence_equilibrium} as $n\to \infty$ and prove that the limit satisfies a particular ODE thanks to Proposition \ref{prop_Poisson_psi}, which will be referred as Steps 1 and 2 in the proofs in the forthcoming Sections \ref{sec:fast} and \ref{sec:equilibrium}.
\item Then, we will identify this limit as the LT of the right hand side of \eqref{convergence_fast} and \eqref{convergence_equilibrium} thanks to a proper use of the Feynman-Kac or Campbell formula. This step will be referred as Step 3 in the proofs.
\end{itemize}
This is to be compared with the approach in \cite{BKMT14, MDT16}, where the authors derive ODEs for the limiting moment generating function and identify a gaussian limiting distribution for the normalized process.

The paper is organized in the following way. Section \ref{main_proofs} is dedicated to the proofs of the main results, with Subsections \ref{sec:fast}, \ref{sec:equilibrium} and \ref{sec:slow} giving the proofs of the convergences in distribution of Theorem \ref{theo_regimes} in fast arrivals and equilibrium cases, and of Theorem \ref{theo_slow_arrival} in the slow arrivals case. The proof of Corollary \ref{rem:marginal} is included  in Subsections \ref{sec:fast}, for the convergence \eqref{remark_conv_distrib_simpler_fast}, and \ref{sec:equilibrium}, for the convergence \eqref{remark_conv_distrib_simpler_equilibirum}. As a concluding remark, we will discuss in Section \ref{sec:remark_compute} some computational aspect for the limiting distributions mentioned in those different regimes in Theorem \ref{theo_regimes} in the particular case when $\alpha$ is a rational number lying in $(0,1)$.

\section{Proofs of Theorems \ref{theo_regimes}, \ref{theo_slow_arrival} and Corollary \ref{rem:marginal}}\label{main_proofs}

\subsection{Preliminary results}\label{sec:preliminary}
We will repeatedly use the following general lemma in the proofs:
\begin{lemm}\label{lemma_convergence}
Let $\left(t\in [0,1]\mapsto A_n(t) \right)_{n\in\nbN}$ be a sequence of continuous functions with values in $\nbR^{{\cal S}\times {\cal S}}$, and let us assume that there exists some continuous function $t\in [0,1]\mapsto A(t)\in \nbR^{{\cal S}\times {\cal S}}$ such that $\int_0^1 || A_n(v)-A(v)|| dv \longrightarrow 0$ as $n\to\infty$ for any matrix norm $||.||$. Let $y\in  \nbR^{{\cal S}\times {\cal S}}$ and $t\in [0,1]\mapsto Y_n(t)\in \nbR^{{\cal S}\times {\cal S}}$ be the solution to the following differential equation
\begin{equation}\label{lemma_EDOn}
\left\{
\begin{array}{rcl}
\frac{d}{dt}Y_n(t) &=& A_n(t) Y_n(t),\quad t\in [0,1],\\
Y_n(0)&=& y\in \nbR^{{\cal S}\times {\cal S}} ,
\end{array}
\right. \quad n\in \nbN .
\end{equation}
Then one has $Y_n(t)\longrightarrow Y(t)$ uniformly in $t\in[0,1]$, as $n\to \infty$, where $t\in [0,1]\mapsto Y(t)\in \nbR^{{\cal S}\times {\cal S}}$ is the solution to the following differential equation
\begin{equation}\label{lemma_EDO}
\left\{
\begin{array}{rcl}
\frac{d}{dt}Y(t) &=& A(t) Y(t),\quad t\in [0,1],\\
Y(0)&=& y .
\end{array}
\right.
\end{equation}
\end{lemm}
\begin{proof}
We first observe that, because of continuity of $t\in [0,1]\mapsto A_n(t)$ and $t\in [0,1]\mapsto A(t)$, \eqref{lemma_EDOn} and \eqref{lemma_EDO} read in integral form
\begin{equation}\label{equ_diff_Y}
Y_n(t)= y + \int_0^t A_n(v) Y_n(v) dv,\quad Y(t)= y + \int_0^t A(v) Y(v) dv
\end{equation}
for all $t\in [0,1]$. Since the norm $||.||$ may be arbitrary, we pick a submultiplicative one on the set of ${\cal S}\times {\cal S}$ matrices. \eqref{equ_diff_Y} implies the following inequality
$$
||Y_n(t)|| \le ||y||+ \int_0^t ||A_n(v)||.|| Y_n(v)|| dv,\quad \forall t\in [0,1] .
$$
Gronwall's lemma thus implies that $ ||Y_n(t)|| \le ||y|| \exp\left( \int_0^t ||A_n(v)|| dv \right)$ for all $t\in [0,1]$. Since by assumption $\int_0^1 || A_n(v)-A(v)|| dv \longrightarrow 0$ as $n\to\infty$, one has that $\left(\int_0^1 ||A_n(v)|| dv\right)_{n\in \N}$ is a bounded sequence. We deduce the following finiteness
\begin{multline*}
M_Y:= \sup_{n\in\N}\sup_{t\in [0,1]}||Y_n(t)|| \le   \sup_{n\in\N}\sup_{t\in [0,1]} ||y|| \exp\left( \int_0^t ||A_n(v)|| dv \right)\\
\le   ||y|| \exp\left( \int_0^1 \sup_{n\in\N} ||A_n(v)|| dv \right) <+\infty .
\end{multline*}
Let us then introduce $M_A:= \sup_{v\in [0,1]}||A(v)||$, which is a finite quantity. Then one obtains that
\begin{multline*}
||Y_n(t)-Y(t)||\le \int_0^t ||A_n(v)-A(v)||. ||Y_n(v)|| dv +  \int_0^t ||A(v)||.||Y_n(v)-Y(v)||dv\\
\le M_Y \int_0^t ||A_n(v)-A(v)|| dv + M_A \int_0^t ||Y_n(v)-Y(v)||dv,\quad \forall t\in[0,1].
\end{multline*}
Gronwall's lemma thus implies that, for all $t\in [0,1]$,
\begin{multline*}
||Y_n(t)-Y(t)||\le M_Y \left[ \int_0^t ||A_n(v)-A(v)|| dv\right].\ e^{M_A t}\\ \le M_Y \left[  \int_0^1 ||A_n(v)-A(v)|| dv \right].\ e^{M_A}\longrightarrow 0 \mbox{ as } n\to \infty.
\end{multline*}
Since the right hand side of the above inequality is independent from $t\in [0,1]$, this proves the uniform convergence result.
\end{proof}
We finish this subsection by stating the differential equation satisfied by the Laplace Transform $\psi^{(n)}(s,t)$ of $\left(Z^{(n)}(t),X^{(n)}_{N^{(n)}_t}\right)$ defined in \eqref{def_LT_n}, which will be the central object studied in Subsections \ref{sec:fast} and \ref{sec:equilibrium}.
Thanks to equation\eqref{Poisson_ODE} with the new parameters $\lambda_n$, $P_n$ instead of $\lambda$ and $P$ (and remembering that $\lambda_n (P_n-I)=\lambda(P-I)$), this reads here
\begin{equation}\label{Poisson_ODEn}
\left\{
\begin{array}{rcl}
\partial_t \psi^{(n)}(s,t) &=& [\lambda (P-I) + \lambda n^\gamma P_n(\tilde{\pi}_n(s,t)-I)]\psi^{(n)}(s,t),\quad t\ge 0,\\
\psi^{(n)}(s,0)&=& I,
\end{array}
\right.
\end{equation}
for all $s=(s_1,...,s_k)\in (-\infty,0]^k$. And, from \eqref{def_pi_Q_tilda}, using the expansion $\prod_{j=1}^k (a_j+1)=1+ \sum_{I\subset \{1,...,k\}} \prod_{\ell \in I} a_\ell$ for all real numbers $a_1,...,a_k$, we have the following expansion which will be useful later on:
\begin{multline}\label{def_pi_n}
\tilde{\pi}_n(s,t)-I=\mbox{diag} \left( \prod_{j=1}^k \left( (e^{s_j x_j}-1)\nbP \left[L_j^{(n)}>t\right]+1\right)-1, \ x=(x_1,...,x_k)\in {\cal S}\right)\\
= \mbox{diag} \left(  \sum_{I\subset \{1,...,k\}} \prod_{\ell \in I}\left[(e^{s_\ell x_\ell}-1)\nbP\left[L_\ell^{(n)}>t\right]\right], \ x=(x_1,...,x_k)\in {\cal S}\right).
\end{multline}

\subsection{Case $\gamma>\alpha$: Fast arriving customers}\label{sec:fast}
We now proceed to show convergence \eqref{convergence_fast} in Theorem \ref{theo_regimes}. In the present case, it is sensible to guess that $Z^{(n)}(t)$ converges towards infinity as $n\to\infty$, hence it is natural to find a normalization such that a convergence towards a proper distribution occurs. We renormalize the queue content by dividing it by $n^{\gamma-\alpha}$, i.e. we are here interested in $\left(Z^{(n)}(t)/n^{\gamma-\alpha},X^{(n)}_{N^{(n)}_t} \right)$, of which Laplace transform is given by $\psi^{(n)}(s/n^{\gamma-\alpha},t)$, $s=(s_1,...,s_k)\in (-\infty,0]^k$. In order to avoid cumbersome notation, we introduce the quantity
$$
\beta:= \frac{1}{1-\alpha}\in (1,+\infty).
$$
We observe then that 
\begin{equation}\label{time_transfo}
t\in[0,1]\mapsto t^\beta \in [0,1]
\end{equation}
is a one to one mapping. Hence, studying the limiting distribution of $\left(Z^{(n)}(t)/n^{\gamma-\alpha},X^{(n)}_{N^{(n)}_t} \right)$ for all $t\in [0,1]$ amounts to study the limiting distribution of
\begin{equation}\label{fast_renormalized_beta}
\left(Z^{(n)}(t^\beta)/n^{\gamma-\alpha},X^{(n)}_{N^{(n)}_{t^\beta }} \right)
\end{equation}
for all $t\in [0,1]$, then changing variable $t:=t^{1/\beta}$. The time transformation \eqref{time_transfo} may at this point look artificial, but this is a key step which will later on enable us to use the convergence result in Lemma \ref{lemma_convergence}.The LT of \eqref{fast_renormalized_beta} is given by
$$\chi^{(n)}(s,t):= \psi^{(n)}(s/n^{\gamma-\alpha},t^\beta),\quad t\in [0,1].$$
From \eqref{Poisson_ODEn}, $\chi^{(n)}(s,t)$ satisfies
\begin{equation}\label{Poisson_ODEn_fast}
\left\{
\begin{array}{rcl}
\partial_t \chi^{(n)}(s,t) &=& \beta t^{\beta-1}[\lambda (P-I) + \lambda n^\gamma P_n(\tilde{\pi}_n(s/n^{\gamma-\alpha},t^\beta)-I)] \chi^{(n)}(s,t),\quad t\in [0,1],\\              
\chi^{(n)}(s,0)&=& I.
\end{array}
\right.
\end{equation}
The starting point for proving \eqref{convergence_fast} is the following: we will set to prove that 
\begin{equation}\label{def_An_fast}
A_n(s,t)=\beta t^{\beta-1}[\lambda (P-I) + \lambda n^\gamma P_n(\tilde{\pi}_n(s/n^{\gamma-\alpha},t^\beta)-I)]
\end{equation}
converges to some limit $A(s,t)$ as $n\to\infty$, use Lemma \ref{lemma_convergence}, then identify the limit $\chi(s,t):=\lim_{n\to\infty}\chi^{(n)}(s,t)$ as the Laplace Transform of a known distribution.\\
{\bf Step 1: Determining $A(s,t)$.} This step is dedicated to finding the limit function $t\in[0,1]\mapsto A(s,t)$ of \eqref{def_An_fast}. Recalling that $\lim_{n\to \infty}P_n=I$, studying the limit of \eqref{def_An_fast} amounts to studying that of $\beta t^{\beta-1}\lambda n^\gamma(\tilde{\pi}_n(s/n^{\gamma-\alpha},t^\beta)-I)$ as $n\to\infty$. In view of \eqref{def_pi_n}, the $x$th diagonal element of this latter term is
\begin{equation}\label{term1_fast}
\beta t^{\beta-1}\lambda n^\gamma \sum_{I\subset \{1,...,k\}} \prod_{\ell \in I}\left[(e^{s_\ell x_\ell/n^{\gamma-\alpha}}-1)\nbP\left[L_\ell^{(n)}>t^\beta\right]\right]
\end{equation}
of which we proceed to find the limit as $n\to \infty$. In order to study its convergence, we are going to isolate the terms in the sum \eqref{term1_fast} for which $\mbox{Card}(I)=1$ and $\mbox{Card}(I)\ge 2$, and show that the former admit a non zero limit and the latter tend to $0$. We thus write \eqref{term1_fast} as
\begin{eqnarray}
&& \beta t^{\beta-1}\lambda n^\gamma \sum_{I\subset \{1,...,k\}} \prod_{\ell \in I}\left[(e^{s_\ell x_\ell/n^{\gamma-\alpha}}-1)\nbP\left[L_\ell^{(n)}>t^\beta\right]\right]=J_n^1(s,t) + J_n^2(s,t),\quad \mbox{where}\nonumber\\
&&J_n^1(s,t)=J_n^1(s,t,x):= \beta t^{\beta-1}\lambda n^\gamma \sum_{\ell =1}^k (e^{s_\ell x_\ell/n^{\gamma-\alpha}}-1)\nbP\left[L_\ell^{(n)}>t^\beta\right],\label{def_J1n}\\
&&J_n^2(s,t)=J_n^2(s,t,x):= \beta t^{\beta-1}\lambda n^\gamma \sum_{\mbox{\tiny Card}(I)\ge 2} \ \prod_{\ell \in I}\left[(e^{s_\ell x_\ell/n^{\gamma-\alpha}}-1)\nbP\left[L_\ell^{(n)}>t^\beta\right]\right].\label{def_J2n}
\end{eqnarray}
Both terms $J_n^1(s,t)$ and $J_n^2(s,t)$ are studied separately. Using that $e^{s_\ell x_\ell/n^{\gamma-\alpha}}-1\sim s_\ell x_\ell/n^{\gamma-\alpha}$ as $n\to\infty$ and 
\begin{equation}\label{basic_fat_tail}
\nbP\left[L_\ell^{(n)}>t^\beta\right]=\nbP\left[L_\ell>nt^\beta\right]\sim \frac{1}{n^\alpha t^{\beta \alpha}}
\end{equation}
when $t>0$, and since $\beta\alpha=\alpha/(1-\alpha)=\beta-1$, we arrive at
$$
J_n^1(s,t)\sim \beta \lambda \sum_{\ell =1}^k t^{\beta-1} n^\gamma \frac{s_\ell x_\ell}{n^{\gamma-\alpha}} \frac{1}{n^\alpha t^{\beta \alpha}}\sim \beta \lambda \sum_{\ell =1}^k s_\ell x_\ell,\quad n\to \infty, 
$$
when $t>0$, and is $0$ when $t=0$.
Next we show that $J_n^2(s,t)$ tends to $0$ by showing that each term on the right hand side of \eqref{def_J2n} tend to $0$. So, if $I\subset \{1,...,k\}$ is such that $I=\{\ell_1,\ell_2\}$, i.e. $\mbox{Card}(I) =2$, then 
\begin{multline}\label{J2n_tends_zero}
\left| \beta t^{\beta-1}\lambda n^\gamma \prod_{\ell \in I}\left[(e^{s_\ell x_\ell/n^{\gamma-\alpha}}-1)\nbP\left[L_\ell^{(n)}>t^\beta\right]\right]\right|= \beta t^{\beta-1}\lambda n^\gamma \left|e^{s_{\ell_1} x_{\ell_1}/n^{\gamma-\alpha}}-1\right|\nbP\left[L_{\ell_1}>nt^\beta\right]\\
. \left|e^{s_{\ell_2} x_{\ell_2}/n^{\gamma-\alpha}}-1\right|\nbP\left[L_{\ell_2}>nt^\beta\right]\le \beta t^{\beta-1}\lambda n^\gamma |s_{\ell_1} x_{\ell_1}|.|s_{\ell_2} x_{\ell_2}| \frac{1}{n^{2(\gamma-\alpha)}}\nbP\left[L_{\ell_1}>nt^\beta\right]\\
= \beta t^{\beta-1}\lambda |s_{\ell_1} x_{\ell_1}|.|s_{\ell_2} x_{\ell_2}| \frac{1}{n^{\gamma-\alpha}} n^\alpha \nbP\left[L_{\ell_1}>nt^\beta\right],
\end{multline}
where we used the inequality $|e^u-1|\le |u|$ for $u\le 0$ and $\nbP\left[L_{\ell_2}>nt^\beta\right]\le 1$. Thanks to \eqref{basic_fat_tail}, the right hand side of \eqref{J2n_tends_zero} thus tends to zero when $t\in (0,1]$. The case $\mbox{Card}(I) >2$ is dealt with similarly. Finally, all terms on the right hand side of \eqref{def_J2n} tend to $0$ as $n\to\infty$, i.e. $\lim_{n\to\infty}J_n^2(s,t)=0$ for all $t\in (0,1]$. When $t=0$ then $J_n^2(s,t)=0$, so that the limit holds for all $t\in [0,1]$.\\
Hence we have that \eqref{term1_fast} tends to $\lim_{n\to\infty}J_n^1(s,t)+\lim_{n\to\infty}J_n^2(s,t)$, i.e. to $ \beta \lambda\sum_{\ell =1}^k s_\ell x_\ell$ when $t\in (0,1]$, and to $0$ when $t=0$. The candidate for the continuous function $A(s,t)$ is then
\begin{equation}\label{candidate_Ast}
t\in [0,1]\mapsto A(s,t):= \beta t^{\beta-1}\lambda (P-I) + \beta \lambda \sum_{\ell =1}^k s_\ell\Delta_\ell
\end{equation}
where we recall from \eqref{Di} that $\Delta_\ell=\mbox{diag} \left[ x_\ell ,\ x=(x_1,\ldots,x_k)\in{\cal S}\right]$. This is where the time transformation \eqref{time_transfo} described previously is important, as without it it would not have been possible to exhibit the limit \eqref{candidate_Ast} for $A_n(s,t)$. Note that the limit when $t=0$ for $A_n(s,t)$ in \eqref{def_An_fast} differs from $A(s,0)=\beta \lambda \sum_{\ell =1}^k s_\ell\Delta_\ell$, as indeed a closer look from the study of the limits of $J_n^1(s,t)$ and $J_n^2(s,t)$ would yield that $\lim_{n\to \infty}A_n(s,0)$ should rather be the $0$ matrix. This is due to the fact that the limit $t\in [0,1]\mapsto A(s,t)$ in Lemma \ref{lemma_convergence} has to be {\it continuous} so that the lemma holds.\\
{\bf Step 2: Determining $\chi(s,t)=\lim_{n\to}\chi_n(s,t)$.} We now need to prove that $\int_0^1 || A_n(s,v)-A(s,v)|| dv \longrightarrow 0$ as $n\to\infty$ in order to apply Lemma \ref{lemma_convergence}. Thanks to \eqref{def_An_fast} and \eqref{candidate_Ast}, and by the definitions \eqref{def_J1n} and \eqref{def_J2n} of $J^1_n(s,t,x)$ and $J^1_n(s,t,x)$, we observe that $A_n(s,t)$ can be decomposed as
\begin{multline*}
A_n(s,t)=A(s,t) + P_n\ \mbox{diag}\left(J_n^1(s,t,x) - \beta \lambda \sum_{\ell =1}^k s_\ell x_\ell,\ x\in {\cal S}\right) \\
+ P_n\ \mbox{diag}\left(J_n^2(s,t,x),\ x\in {\cal S}\right) + (P_n-I)\beta \lambda \sum_{\ell =1}^k s_\ell\Delta_\ell, \quad t\in[0,1].
\end{multline*}
Hence, since $\lim_{n\to \infty}P_n=I$, proving $\lim_{n\to\infty}\int_0^1 || A_n(s,v)-A(s,v)|| dv = 0$ amounts to proving that
\begin{equation}\label{to_prove_limits_J12}
\begin{array}{rcl}
\displaystyle\int_0^1 \left|J_n^1(s,v,x) - \beta \lambda \sum_{\ell =1}^k s_\ell x_\ell\right| dv & \longrightarrow & 0, \ \mbox{and}\\
\displaystyle\int_0^1 |J_n^2(s,v,x)|dv=\int_0^1 J_n^2(s,v,x)dv & \longrightarrow & 0,
\end{array}
\end{equation}
as $n\to\infty$, for each fixed $x\in {\cal S}$. Let us first focus on $\int_0^1 \left|J_n^1(s,v,x) - \beta \lambda \sum_{\ell =1}^k s_\ell x_\ell\right| dv$. We have
\begin{eqnarray}
\int_0^1 \left|J_n^1(s,v,x) - \beta \lambda \sum_{\ell =1}^k s_\ell x_\ell\right| dv&\le & \sum_{\ell=1}^k  (I^1_n(\ell) + I^2_n(\ell)), \mbox{ where, for all } \ell =1,...,k, \label{term3_fast}\\
I^1_n(\ell) &:=& \int_0^1  \lambda \beta v^{\beta-1}\left|  n^\gamma(e^{s_\ell x_\ell/n^{\gamma-\alpha}}-1)-   n^\alpha s_\ell x_\ell\right| \nbP\left[L_\ell^{(n)}>v^\beta\right]dv \nonumber\\
I^2_n(\ell) &:=& |s_\ell x_\ell| \int_0^1  \lambda  \left|\beta v^{\beta-1} n^\alpha   \nbP\left[L_\ell^{(n)}>v^\beta\right]-\beta \right|dv .\nonumber
\end{eqnarray}
Expanding the exponential function, one has that $|e^{s_\ell x_\ell/n^{\gamma-\alpha}}-1-s_\ell x_\ell/n^{\gamma-\alpha}|\le M_\ell/n^{2(\gamma-\alpha)}$ where $M_\ell>0$ only depends on $s_\ell$ and $x_\ell$. Thus, one deduces the following upper bounds for $I^1_n(\ell)$, $\ell =1,...,k$:
\begin{eqnarray}
I^1_n(\ell) &=& \int_0^1   \lambda \beta v^{\beta-1}\left|  n^\gamma(e^{s_\ell x_\ell/n^{\gamma-\alpha}}-1)-   n^\alpha s_\ell x_\ell\right| \nbP\left[L_\ell> n v^\beta\right]dv\nonumber\\
&\le & n^\gamma \frac{M_\ell}{n^{2(\gamma-\alpha)}} \lambda \int_0^1 \beta v^{\beta-1} \nbP\left[L_\ell> n v^\beta\right]dv = \frac{M_\ell}{n^{\gamma-\alpha}} \beta \lambda \int_0^1 n^\alpha v^{\beta-1} \nbP\left[L_\ell> n v^\beta\right]dv \nonumber\\
 &=& \frac{M_\ell}{n^{\gamma-\alpha}} \beta \lambda \int_0^1 (nv^{\beta})^\alpha \ \nbP\left[L_\ell> n v^\beta\right]dv,\label{term4_fast}
\end{eqnarray}
the last equality holding because $\beta-1=\beta \alpha$ implies that the integrand verifies $n^\alpha v^{\beta-1}  =(nv^{\beta})^\alpha$. A consequence of the fact that $L_\ell$ is fat tailed with index $\alpha$ is that $\sup_{u\ge 0}u^\alpha \nbP(L_\ell>u)<+\infty$, from which one deduces immediately that
\begin{equation}\label{bound_sup_L_l}
\sup_{j\in\nbN,\ v \in [0,1]} (jv^{\beta})^\alpha \ \nbP\left[L_\ell> j v^\beta\right]<+\infty
\end{equation}
(note that those two latter suprema are in fact equal). One then gets from \eqref{term4_fast} that
\begin{equation}\label{term5_fast}
I^1_n(\ell)\le \frac{M_\ell}{n^{\gamma-\alpha}} \beta \lambda \left[\sup_{j\in\nbN,\ v \in [0,1]} (jv^{\beta})^\alpha \ \nbP\left[L_\ell> j v^\beta\right]\right]\longrightarrow 0,\quad n\to\infty .
\end{equation}
We now turn to $I^2_n(\ell)$, $\ell =1,...,k$. Using again $\beta-1=\beta \alpha$, one may write in the integrand of $I^2_n(\ell)$ that $v^{\beta-1} n^\alpha =(nv^{\beta})^\alpha$: hence
$$ I^2_n(\ell) = |s_\ell x_\ell| \int_0^1  \lambda  \left|\beta (nv^{\beta})^\alpha   \nbP\left[L_\ell>nv^\beta\right]-\beta \right|dv .$$
Since $L_\ell$ is fat tailed with index $\alpha$, estimates similar to the ones leading to the upper bound \eqref{term5_fast} for $I_n^1(\ell)$ yield that
$$
\sup_{n\in\nbN,\ v \in [0,1]}\left|\beta (nv^{\beta})^\alpha   \nbP\left[L_\ell>nv^\beta\right]-\beta \right|<+\infty .
$$
Furthermore, again because $L_\ell$ is fat tailed, one has $\nbP\left[L_\ell>nv^\beta\right] \sim 1/(nv^{\beta})^\alpha$ as $n\to \infty$ when $v>0$. Hence $\left|\beta (nv^{\beta})^\alpha   \nbP\left[L_\ell>nv^\beta\right]-\beta \right|\longrightarrow 0 $ as $n\to \infty$ when $v\in (0,1]$, and is equal to $\beta$ when $v=0$. The dominated convergence theorem thus implies that
\begin{equation}\label{term6_fast}
I^2_n(\ell)\longrightarrow 0, \quad n\to\infty .
\end{equation}
Gathering \eqref{term3_fast}, \eqref{term5_fast} and \eqref{term6_fast}, we thus deduce finally that $\int_0^1 \left|J_n^1(s,v,x) - \beta \lambda \sum_{\ell =1}^k s_\ell x_\ell\right| dv$ tends to $0$ as $n\to\infty$ for each $x\in{\cal S}$. \\
We now prove that $\int_0^1 J_n^2(s,v,x)dv\longrightarrow 0$ as $n\to \infty$.
In view of the definition \eqref{def_J2n}, it suffices to prove that 
\begin{equation}\label{to_prove_J2}
\int_0^1\beta v^{\beta-1}\lambda n^\gamma  \ \prod_{\ell \in I}\left[|e^{s_\ell x_\ell/n^{\gamma-\alpha}}-1|\ \nbP\left[L_\ell^{(n)}>v^\beta\right]\right]dv
\end{equation}
tends to $0$ as $n\to\infty$ for $I\subset \{1,...,k\}$ such that $\mbox{ Card}(I)\ge 2$. Let us prove the convergence for $\mbox{ Card}(I)= 2$, i.e. for $I=\{\ell_1,\ell_2\}$ for some $\ell_1\neq\ell_2$ in $1,...,k$, the case $\mbox{ Card}(I)> 2$ being dealt with similarly. By the basic inequality $|e^u-1|\le |u|$ for $u\le 0$ we deduce that $|e^{s_{\ell_i} x_{\ell_i}/n^{\gamma-\alpha}}-1|\le |s_{\ell_i} x_{\ell_i}|/n^{\gamma-\alpha}$, $i=1,2$. Since $\nbP\left[L_{\ell_1}^{(n)}>v^\beta\right]\le 1$ for all $v\in [0,1]$, we then deduce that \eqref{to_prove_J2} is upper bounded by
$$\frac{|s_{\ell_1} x_{\ell_1}|}{n^{\gamma-\alpha}} |s_{\ell_2} x_{\ell_2}|
\int_0^1 \beta v^{\beta-1}\lambda n^\alpha \nbP\left[L_{\ell_2}^{(n)}>v^\beta\right] dv.
$$
As $v^{\beta-1} n^\alpha =(nv^{\beta})^\alpha$, and thanks to \eqref{bound_sup_L_l}, the latter quantity is in turn written then bounded as follows
\begin{multline*}
\frac{|s_{\ell_1} x_{\ell_1}|}{n^{\gamma-\alpha}} |s_{\ell_2} x_{\ell_2}|
\int_0^1 \beta \lambda  (nv^{\beta})^\alpha \nbP\left[L_{\ell_2}^{(n)}>v^\beta\right] dv
= \frac{|s_{\ell_1} x_{\ell_1}|}{n^{\gamma-\alpha}} |s_{\ell_2} x_{\ell_2}|
\int_0^1 \beta \lambda  (nv^{\beta})^\alpha \nbP\left[L_{\ell_2}>nv^\beta\right] dv\\
\le \frac{|s_{\ell_1} x_{\ell_1}|}{n^{\gamma-\alpha}} |s_{\ell_2} x_{\ell_2}| \beta \lambda \left[\sup_{j\in\nbN,\ v \in [0,1]} (jv^{\beta})^\alpha \ \nbP\left[L_{\ell_2}> j v^\beta\right] \right] \longrightarrow 0,\quad n\to\infty, 
\end{multline*}
proving that \eqref{to_prove_J2} tends to $0$ as $n\to\infty$ when $I=\{\ell_1,\ell_2\}$.\\
Hence we just proved \eqref{to_prove_limits_J12}, which implies $\int_0^1 || A_n(s,v)-A(s,v)|| dv \longrightarrow 0$. We may then use Lemma \ref{lemma_convergence} to deduce that $\chi^{(n)}(s,t)$ convgerges to $\chi(s,t)$ which satisfies
\begin{equation}\label{Poisson_ODE_fast}
\left\{
\begin{array}{rcl}
\partial_t \chi(s,t) &=&   A(s,t)\chi(s,t)=\left[\beta t^{\beta-1}\lambda (P-I) + \beta \lambda \sum_{\ell =1}^k s_\ell\Delta_\ell \right]\chi(s,t) ,\quad t\in [0,1],\\
\chi(s,0)&=& I.
\end{array}
\right.
\end{equation}
{\bf Step 3: Identifying the limit in distribution.} Let us note that \eqref{Poisson_ODE_fast} does not admit an explicit expression. However, since we purposely chose $s=(s_1,...,s_k)$ with $s_j\le 0$, $j=1,...,k$, one has that $\sum_{j=1}^k s_j \Delta_j = \sum_{j=1}^k s_j\ \mbox{diag}(x_j,\ x\in {\cal S})$ is a diagonal matrix with non positive entries. Let $\Delta_\pi:=\mbox{\normalfont diag}(\pi(x),\ x\in {\cal S})$ and let us introduce the matrix $P^{(r)}$ defined by $P^{(r)}=\Delta_\pi^{-1}P' \Delta_\pi\iff P=\Delta_\pi^{-1}P^{(r)'} \Delta_\pi$. It is standard that $P^{(r)}$ is the transition matrix of the reversed version of the stationary Markov chain $\{X_i,\ i\in\nbN \}$ with distribution $\pi$, and that $\beta t^{\beta-1}\lambda (P^{(r)}-I)$ is the infinitesimal generator matrix of an inhomogeneous Markov process
\begin{equation}\label{def_U}
\{U(t)=(U_j(t))_{j=1,...,k}\in {\cal S},\ t\in [0,1]\} 
\end{equation}
with values in ${\cal S}$, and initial distribution $U(0)\sim\pi$. 
In fact, it turns out that the conditional distribution of $U(t)$ given $U(0)$ is given by $\left[\nbP(U(t)=y|\ U(0)=x)\right]_{(x,y)\in {\cal S}}=\exp(t^\beta \lambda(P^{(r)}-I))$, which results in $U(t)\sim\pi$ for all $t\in[0,1]$, i.e. that $\{U(t),\ t\in [0,1]\}$ is stationary. Since $\sum_{j=1}^k s_j \Delta_j$ is diagonal, one checks easily that  $A(s,t)=\Delta_\pi^{-1} \left[\beta t^{\beta-1}\lambda ({P^{(r)}}'-I)+\sum_{j=1}^k s_j \Delta_j\right]\Delta_\pi$ and that $Y(t)=Y(s,t):=\Delta_\pi^{-1} \chi(s,t)'\Delta_\pi $ satisfies the differential equation
$$
\left\{
\begin{array}{rcl}
\partial_t Y(t) &=&  Y(t)\left[\beta t^{\beta-1}\lambda ({P^{(r)}}-I) + \beta \lambda \sum_{\ell =1}^k s_\ell\Delta_\ell \right]  ,\quad t\in [0,1],\\
Y(0)&=& I.
\end{array}
\right.
$$
The Feynman-Kac formula ensures that one has the representation
$$Y(t)=Y(s,t)=\left[ \nbE\left[\left. {\bf 1}_{[U(t)=y]}\exp\left(\sum_{j=1}^k s_j \beta \lambda\int_0^t U_j(v) dv  \right)\right| U(0)=x\right]\right]_{(x,y)\in {\cal S}^2},\quad \forall t\in [0,1],$$
see \cite[Chapter III, 19, p.272]{Rogers_Williams00} for the general theorem on this formula, or \cite[Section 5, Expression (5.2) and differential equation (5.3)]{BH96} for the particular case of a finite Markov chain, adapted here to an inhomogeneous Markov process. Also, the reversed process $\{U(1-t),\ t\in [0,1]\}$ admits ${\Delta_\pi}^{-1}\beta (1-t)^{\beta-1}\lambda ({P^{(r)}}'-I)\Delta_\pi= \beta (1-t)^{\beta-1}\lambda (P-I)$ as infinitesimal generator matrix, which is the generator of the process $\{{\cal X}^\alpha (t)=({\cal X}_1^\alpha (t),...,{\cal X}_k^\alpha (t))\in {\cal S},\ t\in[0,1] \}$ introduced in the statement of Theorem \ref{theo_regimes}, so that $\{{\cal X}^\alpha (t),\ t\in[0,1] \}\stackrel{\cal D}{=}\{U(1-t),\ t\in [0,1]\}$ pathwise. Hence, one obtains for all $x$ and $y$ in ${\cal S}$ that
\begin{eqnarray}
&&\nbE\left[\left. {\bf 1}_{[U(t)=y]}\exp\left(\sum_{j=1}^k s_j \beta \lambda\int_0^t U_j(v) dv  \right)\right| U(0)=x\right]\nonumber\\
&=& \nbE\left[\left. {\bf 1}_{[{\cal X}^\alpha (1-t)=y]}\exp\left(\sum_{j=1}^k s_j \beta \lambda\int_{1-t}^1 {\cal X}^\alpha_j (v) dv  \right)\right|{\cal X}^\alpha (1)=x\right]\nonumber\\
&=& \nbE\left[\left. {\bf 1}_{[{\cal X}^\alpha (1)=x]}\exp\left(\sum_{j=1}^k s_j \beta \lambda\int_{1-t}^1 {\cal X}^\alpha_j (v) dv  \right)\right|    {\cal X}^\alpha (1-t)=y\right] \frac{\pi(y)}{\pi(x)},\label{relation_U_chi}
\end{eqnarray}
the last line coming from the fact that $U(0)$, $U(t)$, $ {\cal X}^\alpha (1-t)$ and ${\cal X}^\alpha (1)$ all have same distribution $\pi$. Switching the role of $x$ and $y$ in the above results in the following relationship:
\begin{eqnarray*}
&&\left[ \nbE\left[\left. {\bf 1}_{[{\cal X}^\alpha (1)=y]}\exp\left(\sum_{j=1}^k s_j \beta \lambda\int_{1-t}^1 {\cal X}^\alpha_j (v) dv  \right)\right|    {\cal X}^\alpha (1-t)=x\right] \right]_{(x,y)\in {\cal S}^2}\\
&=& \left[ \nbE\left[\left. {\bf 1}_{[U(t)=x]}\exp\left(\sum_{j=1}^k s_j \beta \lambda\int_0^t U_j(v) dv  \right)\right| U(0)=y\right] \frac{\pi(y)}{\pi(x)}\right]_{(x,y)\in {\cal S}^2}\\
&=& \Delta_\pi^{-1} Y(t)'\Delta_\pi=\chi(s,t).
\end{eqnarray*}
Since we just proved that $\chi^{(n)}(s,t):= \psi^{(n)}(s/n^{\gamma-\alpha},t^\beta)$ converges as $n\to \infty$ towards $\chi(s,t)$, expressed above, for all $s=(s_1,...,s_k)\in (-\infty, 0]^k$, and identifying Laplace transforms, we obtained in conclusion that 
\begin{equation}\label{conv_fast_final}
{\cal D}\left(\left.  \left(Z^{(n)}(t^\beta)/n^{\gamma-\alpha},X^{(n)}_{N^{(n)}_{t^\beta }} \right)\right|\ X^{(n)}_0 \right)\longrightarrow {\cal D}\left(
\left. 
\left( \beta\lambda \int_{1-t}^{1} {\cal X}^\alpha(v)\ dv,\ {\cal X}^\alpha (1)\right)
\right|\ {\cal X}^\alpha(1-t)
\right)
\end{equation}
as $n\to\infty$ for all $t\in[0,1]$. Changing $t$ into $t^{1/\beta}$ yields \eqref{convergence_fast}.\\
{\bf Proof of the convergence \eqref{remark_conv_distrib_simpler_fast} in Corollary \ref{rem:marginal}. }With the previous definitions of processes $\{U(t),\ t\in [0,1]\}$ in \eqref{def_U} and $\{{\cal X}^\alpha (t),\ t\in[0,1] \}$, \eqref{relation_U_chi} implies the following matrix equality
\begin{multline}\label{relation_U_chi_matrix}
\left[ \nbE\left[\left. {\bf 1}_{[U(t)=y]}\exp\left(\sum_{j=1}^k s_j \beta \lambda\int_0^t U_j(v) dv  \right)\right| U(0)=x\right] \right]_{(x,y)\in {\cal S}^2}\\
= \left[ \nbE\left[\left. {\bf 1}_{[{\cal X}^\alpha (1)=x]}\exp\left(\sum_{j=1}^k s_j \beta \lambda\int_{1-t}^1 {\cal X}^\alpha_j (v) dv  \right)\right|    {\cal X}^\alpha (1-t)=y\right] \frac{\pi(y)}{\pi(x)}\right]_{(x,y)\in {\cal S}^2}
\end{multline}
Left-multiplying and right-multiplying \eqref{relation_U_chi_matrix} respectively by the row vector $(\pi(x))_{x\in {\cal S}}$ and the column vector ${\bf 1}$ results in in the following equality of LT
\begin{multline}\label{LT_reverse_fast}
\nbE\left[ \exp\left(\sum_{j=1}^k s_j \beta \lambda\int_0^t U_j(v) dv  \right) \right]\\
= \nbE\left[\exp\left(\sum_{j=1}^k s_j \beta \lambda\int_{1-t}^1 {\cal X}^\alpha_j (v) dv  \right)\right],\quad s=(s_1,...,s_k)\in (-\infty, 0]^k,
\end{multline}
which, combined with \eqref{conv_fast_final}, yields the convergence $\frac{Z^{(n)}\left(t^\beta\right)}{n^{\gamma-\alpha}}\stackrel{\cal D}{\longrightarrow} \beta\lambda \int_{0}^{t} U(v)\ dv$ as $n\to \infty$. Changing $t$ into $t^{1/\beta}$ and performing the change of variable $ v:=v^{1/\beta}=v^{1-\alpha} $, we obtain
\begin{equation}\label{conv_distribution_reversed_fast}
\frac{Z^{(n)}\left(t\right)}{n^{\gamma-\alpha}}\stackrel{\cal D}{\longrightarrow} \lambda \int_{0}^{t} \frac{U(v^{1-\alpha})}{v^\alpha}\ dv,\ n\to\infty, \quad t\in [0,1].
\end{equation}
Since the ${\cal S}$ valued Markov process $\{U(t),\ t\in [0,1]\} $ admits $\beta t^{\beta-1}\lambda (P^{(r)}-I)$ as the infinitesimal generator matrix, the time changed Markov process $\{{\cal Y}(t):=U(t^{1-\alpha})=U(t^{1/\beta}),\ t\in [0,1]\}  $ admits \eqref{generator_Y} as generator, so that \eqref{remark_conv_distrib_simpler_fast} follows from \eqref{conv_distribution_reversed_fast}.

\subsection{Equilibrium case $\gamma=\alpha$}\label{sec:equilibrium}
We now proceed to show convergence \eqref{convergence_equilibrium} in Theorem \ref{theo_regimes}. Intuitively, we are in the critical case where customers should arrive just fast enough such that the queue at time $t$ converges as $n\to\infty$. We are here interested in the behaviour of ${\cal D}\left(
\left.  \left(Z^{(n)}(t),X^{(n)}_{N^{(n)}_t} \right) \right|\ X^{(n)}_0 \right)$ as $n\to\infty$ when $t\in[0,1]$ is fixed. As in Section \ref{sec:fast}, we first consider $t^\beta$ instead of $t$ and let 
$$\chi^{(n)}(s,t):= \psi^{(n)}(s,t^\beta),\quad t\in [0,1],$$
the corresponding Laplace transform, where $s=(s_1,...,s_k)\in (-\infty,0]^k$. $t\in [0,1]\mapsto \chi^{(n)}(s,t)$ then satisfies, thanks to \eqref{Poisson_ODEn}, the following differential equation
\begin{equation}\label{Poisson_ODEn_equilibrium}
\left\{
\begin{array}{rcl}
\partial_t \chi^{(n)}(s,t) &=& \beta t^{\beta-1}[\lambda (P-I) + \lambda n^\gamma P_n(\tilde{\pi}_n(s,t^\beta)-I)] \chi^{(n)}(s,t),\quad t\in [0,1],\\
\chi^{(n)}(s,0)&=& I.
\end{array}
\right.
\end{equation}
The present case has the same roadmap as Subsection \ref{sec:fast}: we will study the behaviour as $n\to\infty$ of $\lambda n^\gamma(\tilde{\pi}_n(s,t^\beta))-I)$ in order to obtain a limit as $n\to\infty$ of
\begin{equation}\label{def_An_equilibirum}
A_n(s,t)=\beta t^{\beta-1}[\lambda (P-I) + \lambda n^\gamma P_n(\tilde{\pi}_n(s,t^\beta))-I)]
\end{equation}
then getting a limiting matrix differential equation with solution $\chi(s,t)=\lim_{n\to\infty} \chi^{(n)}(s,t)$. Then we will identify $\chi(s,t)$ as the Laplace transform of a (conditional) distribution, yielding \eqref{convergence_equilibrium}.\\
{\bf Step 1: Determining $A(s,t)=\lim_{n\to \infty }A_n(s,t)$.} We recall that the $(x,x)$th diagonal element of $\lambda n^\gamma(\tilde{\pi}_n(s,t^\beta))-I)$ is (from \eqref{def_pi_n})
$\sum_{I\subset \{1,...,k\}} \prod_{\ell \in I}\left[(e^{s_\ell x_\ell}-1)\nbP\left[L_\ell>nt\right]\right]$, which we decompose as in Section \ref{sec:fast} as $K_n^1(s,t)+K_n^2(s,t)$ with
\begin{eqnarray}
 K_n^1(s,t)&=&K_n^1(s,t,x):= \beta t^{\beta-1}\lambda n^\gamma \sum_{\ell =1}^k (e^{s_\ell x_\ell}-1)\nbP\left[L_\ell >n t^\beta\right],\label{def_K1n}\\
K_n^2(s,t)&=&K_n^2(s,t,x):= \beta t^{\beta-1}\lambda n^\gamma \sum_{\mbox{\tiny Card}(I)\ge 2} \ \prod_{\ell \in I}\left[(e^{s_\ell x_\ell}-1)\nbP\left[L_\ell >nt^\beta\right]\right].\label{def_K2n}
\end{eqnarray}
The important point here is that, throughout this subsection, we have $\gamma=\alpha$ in the expressions \eqref{def_An_equilibirum}, \eqref{def_K1n} and \eqref{def_K2n}, which will impact on the convergences and limiting results we are going to prove. Using that $\nbP\left[L_\ell >n t^\beta\right]\sim \frac{1}{n^\alpha}\frac{1}{t^{\alpha\beta}}$, $n\to\infty$, when $t>0$, and since $\alpha\beta=\beta-1$, and $\gamma=\alpha$, one here finds that
$$
K_n^1(s,t)=K_n^1(s,t,x)\longrightarrow \left\{
\begin{array}{cl}
\beta\lambda \sum_{\ell =1}^k (e^{s_\ell x_\ell}-1), & t>0,\\
0, & t=0,
\end{array}\right.
 \quad n\to \infty .
$$
As to $K_n^2(s,t)$, one proves easily that it tends to $0$ as $n\to\infty$ for all $t\in [0,1]$, as the sum in \eqref{def_K2n} is over $\mbox{Card}(I)\ge 2$, and using the fat tailed property of the service times. The candidate for the continuous function is thus
\begin{equation}\label{candidate_Ast_equilibrium}
t\in [0,1]\mapsto A(s,t):= \beta t^{\beta-1}\lambda (P-I) + \beta \lambda \sum_{\ell =1}^k \mbox{diag}(e^{s_\ell x_\ell}-1,\ x=(x_1,...,x_k)\in {\cal S}).
\end{equation}
{\bf Step 2: Determining $\chi(s,t)=\lim_{n\to}\chi_n(s,t)$.} We now wish to apply Lemma \ref{lemma_convergence} and prove that $\int_0^1 || A_n(s,v)-A(s,v)|| dv \longrightarrow 0$ where $A_n(s,t)$ and $A(s,t)$ are defined in \eqref{def_An_equilibirum} and \eqref{candidate_Ast_equilibrium}. The method is very similar as to proving \eqref{to_prove_limits_J12} in Step 2 of Section \ref{sec:fast}, as this is equivalent to proving for all $x\in {\cal S}$ that
\begin{eqnarray}
\int_0^1 \left| K_n^1(s,v,x)- \beta\lambda \sum_{\ell =1}^k (e^{s_\ell x_\ell}-1)\right| dv  &\longrightarrow & 0,\label{limits_K1n}\\
\int_0^1 K_n^2(s,v,x)dv  &\longrightarrow & 0 \label{limits_K2n}
\end{eqnarray}
as $n\to\infty$. In view of the expression \eqref{def_K2n} of $K^2_n(s,t)$, \eqref{limits_K2n} is proved the same way as for proving that $\lim_{n\to\infty}\int_0^1 J_n^2(s,v,x)dv=0$ in Step 2 of Section \ref{sec:fast}. More precisely, it suffices from \eqref{def_K2n} to prove that
\begin{equation}\label{proof_limit_K2n}
\lim_{n\to \infty} \int_0^1  \beta t^{\beta-1}\lambda n^\alpha \prod_{\ell \in I} \nbP\left[L_\ell >nt^\beta\right] dt =0
\end{equation}
for all $I\subset \{1,...,k\}$, $\mbox{Card}(I)\ge 2$. We prove it for $I=\{\ell_1,\ell_2\}$, $\ell_1\neq \ell_2$, the proof for $\mbox{Card}(I)> 2$ being very similar. The trick is again to use that $v^{\beta-1} n^\alpha =(nv^{\beta})^\alpha$ as well as the previously established upper bound \eqref{bound_sup_L_l}, resulting in
\begin{multline*}
\int_0^1  \beta t^{\beta-1}\lambda n^\alpha  \nbP\left[L_{\ell_1} >nt^\beta\right] \nbP\left[L_{\ell_2} >nt^\beta\right]dt \\
\le \beta\lambda \left[\sup_{j\in\nbN,\ v \in [0,1]} (jv^{\beta})^\alpha \ \nbP\left[L_{\ell_1}> j v^\beta\right] \right]\ \int_0^1 \nbP\left[L_{\ell_2} >nt^\beta\right]dt,
\end{multline*}
which converges to zero as $n\to\infty$ by the dominated convergence theorem, proving \eqref{proof_limit_K2n} when $\mbox{Card}(I)= 2$.
As to \eqref{limits_K1n}, this is proved, in view of the expression \eqref{def_K1n} of $K^1_n(s,t)$, by showing that $\int_0^1  \lambda  \left|\beta v^{\beta-1} n^\alpha   \nbP\left[L_\ell>nv^\beta\right]-\beta \right|dv$ tends to $0$ as $n\to\infty$ for all $\ell=1,...,k$, as again we have that $\gamma=\alpha$; However,
this was already proved in Step 2 of Section \ref{sec:fast} when proving that $\lim_{n\to\infty}I_n^2(\ell)=0$, $\ell=1,...,k$, see the arguments leading to the convergence \eqref{term6_fast}. All in all, one has the convergence $\int_0^1 || A_n(s,v)-A(s,v)|| dv \longrightarrow 0$, and Lemma \ref{lemma_convergence} is applicable so that $\chi^{(n)}(s,t)$ convgerges to $\chi(s,t)$ which satisfies
\begin{equation}\label{Poisson_ODE_equilibrium}
\left\{
\begin{array}{rcl}
\partial_t \chi(s,t) &=&   A(s,t)\chi(s,t)=\left[\beta t^{\beta-1}\lambda (P-I) \right.\\
&+&\left. \beta \lambda \sum_{\ell =1}^k \mbox{diag}(e^{s_\ell x_\ell}-1,\ x=(x_1,...,x_k)\in {\cal S}) \right]\chi(s,t) ,\quad t\in [0,1],\\
\chi(s,0)&=& I.
\end{array}
\right.
\end{equation}
{\bf Step 3: Identifying the limit in distribution.} With the same notation as in Step 3 of Section \ref{sec:fast} for process $\{{\cal X}^\alpha (t)=({\cal X}_1^\beta (t),...,{\cal X}_k^\beta (t))\in {\cal S},\ t\in[0,1] \}$, one finds this time that
\begin{equation}\label{equilibrium_chi}
\chi(s,t)= \left[ \nbE\left[\left. {\bf 1}_{[{\cal X}^\alpha (1)=y]}\exp\left(\sum_{j=1}^k \beta \lambda\int_{1-t}^1 \left(e^{s_j{\cal X}^\alpha_j (v)}-1\right) dv  \right)\right|    {\cal X}^\alpha (1-t)=x\right] \right]_{(x,y)\in {\cal S}^2}
\end{equation}
for all $s=(s_1,...,s_k)\in (-\infty,0]^k$. We recall the Campbell formula which states that for all measurable function $f:t\in [0,+\infty) \mapsto f(t)\in \nbR$ such that $\int_0^\infty (e^{f(v)}-1) \xi \ dv$ is finite for some $\xi>0$ then one has the identity
$$
\exp \left(\int_0^\infty \left(e^{f(v)}-1\right) \xi \ dv \right) =\nbE \left[ \exp\left( \int_0^\infty f(v)\ \nu(dv)\right)\right],
$$
where $\{ \nu(x),\ x\ge 0\}$ is a Poisson process with intensity $\xi$, see \cite[Section 3.2]{K93}. Conditioning on $\{{\cal X}^\alpha(v),\ v\in [0,1] \}$, this results in \eqref{equilibrium_chi} being written as
$$
\chi(s,t)= \left[ \nbE\left[\left. {\bf 1}_{[{\cal X}^\alpha (1)=y]}\exp\left(\sum_{j=1}^k  s_j\int_{1-t}^{1} {\cal X}^\alpha_j(v)\ \nu_j^\alpha(dv) \right)\right|    {\cal X}^\alpha (1-t)=x\right] \right]_{(x,y)\in {\cal S}^2}
$$
where $\{\nu_j^\alpha(t),\ t\ge 0 \}$, $j=1,...,k$, are independent Poisson processes with intensities $\beta\lambda=\lambda/(1-\alpha)$, and independent from $\{{\cal X}^\alpha(t),\ t\in[0,1] \}$. Identifying Laplace Transforms, we obtain in conclusion that
\begin{multline}\label{conv_equilibrium_final}
{\cal D}\left(\left.  \left(Z^{(n)}(t^\beta),X^{(n)}_{N^{(n)}_{t^\beta }} \right)\right|\ X^{(n)}_0 \right)\\
\longrightarrow {\cal D}\left(
\left. 
\left( \left( \int_{1-t}^{1} {\cal X}^\alpha_j(v)\ d\nu_j^\alpha(v)\right)_{j=1,...,k},\ {\cal X}^\alpha (1)\right)
\right|\ {\cal X}^\alpha(1-t)
\right)
\end{multline}
as $n\to\infty$ for all $t\in [0,1]$. Changing $t$ into $t^{1/\beta}$ completes the proof of \eqref{convergence_equilibrium}.\\
{\bf Proof of the convergence \eqref{remark_conv_distrib_simpler_equilibirum} in Corollary \ref{rem:marginal}. }This follows the same pattern as the proof of \eqref{remark_conv_distrib_simpler_fast}, to which we refer here. More precisely, one verifies this time that, from \eqref{equilibrium_chi}, the analog of \eqref{LT_reverse_fast} in the Fast arrival case is here
\begin{multline}\label{LT_reverse_equilibrium}
\nbE\left[ \exp\left(\sum_{j=1}^k s_j\int_0^t U_j(v)\ \nu_j^\alpha (dv)  \right) \right]\\
= \nbE\left[\exp\left(\sum_{j=1}^k s_j\int_{1-t}^1  {\cal X}^\alpha_j(v)\ d\nu_j^\alpha(v)  \right)\right],\quad s=(s_1,...,s_k)\in (-\infty, 0]^k,
\end{multline}
which, combined with \eqref{conv_equilibrium_final}, yields the convergence $Z^{(n)}\left(t^\beta\right)\stackrel{\cal D}{\longrightarrow} \left( \int_{0}^{t} U_j(v)\ \nu_j^\alpha (dv)\right)_{j=1,...k}$ as $n\to \infty$. Changing $t$ into $t^{1/\beta}$ and performing the change of variable $ v:=v^{1/\beta}=v^{1-\alpha} $, we obtain
\begin{equation}\label{conv_distribution_reversed_equilibrium}
Z^{(n)}\left(t\right)\stackrel{\cal D}{\longrightarrow} \left( \int_{0}^{t} U_j(v^{1-\alpha})\ \tilde{\nu}_j^\alpha(dv)\right)_{j=1,...,k},\ n\to\infty, \quad t\in [0,1],
\end{equation}
where $\{\tilde{\nu}_j^\alpha(t),\ t\in [0,1] \}$, $j=1,...,k$, are the inhomogeneous independent Poisson processes given by $\tilde{\nu}_j^\alpha(t)=\nu_j^\alpha(t^{1/\beta})=\nu_j^\alpha(t^{1-\alpha})$, i.e. Poisson processes with non constant intensity $\lambda t^{-\alpha}$. Arguing, as in the Fast arrival case, the time changed Markov process $\{{\cal Y}(t):=U(t^{1-\alpha})=U(t^{1/\beta}),\ t\in [0,1]\}  $ admits \eqref{generator_Y} as generator, hence \eqref{remark_conv_distrib_simpler_equilibirum} follows from \eqref{conv_distribution_reversed_equilibrium}.

\subsection{Proof of Theorem \ref{theo_slow_arrival}: Slow arriving customers.}\label{sec:slow}
We now consider the case $\gamma<\alpha$. Section 4.2 of \cite{RW18} provide the first two joint moments of the $Z_j^{(n)}(t)$, $j=1,...,k$, with a particular discount factor $a \ge 0$ (recall the notation \eqref{discounted_Z(t)} in Section \ref{sec:model} for the discounted counterpart of the queueing process \eqref{def_Z_t}). Recalling that the rescaling implies that $\lambda_n(P_n-I)=\lambda(P-I)$, we get from \cite[Theorems 14 and 15 with $a=0$ discount factor]{RW18}, that those moments are given by
\begin{eqnarray}
M_j^{(n)}(t) &=& \lambda_n e^{\lambda t (P-I)}\int_0^t \PP \left(L_j^{(n)}>v\right) e^{-\lambda v (P-I)} \Delta_j P_n e^{\lambda v (P-I)} dv, \label{m1_sec:slow}\\
M_{jj}^{(n)}(t) &=&  \lambda_n e^{\lambda t (P-I)}\int_0^t \PP\left(L_j^{(n)}>v\right) e^{-\lambda v (P-I)}
\Delta_j^2 P_n e^{\lambda v (P-I)}\nonumber\\
& + & 2 \PP\left(L_j^{(n)}>v\right) \Delta_j P_n M_{j}^{(n)}(v)dv,\label{m2_jj_sec:slow}\\
M_{jj'}^{(n)}(t) &=&  \lambda_n e^{\lambda t (P-I)}\int_0^t \PP\left(L_j^{(n)}>v\right)\PP\left(L_{j'}^{(n)}>v\right) e^{-\lambda v (P-I)}
\Delta_j\Delta_{j'} P_n e^{\lambda v (P-I)} \nonumber\\
& + & \PP\left(L_j^{(n)}>v\right) \Delta_j P_n M_{j'}^{(n)}(v)
 +  \PP\left(L_{j'}^{(n)}>v\right) \Delta_{j'} P_n M_{j}^{(n)}(v) dv, \label{m2_jj'_sec:slow}
\end{eqnarray}
for all $t\ge 0$ and $j\neq j'$, $j$ and $j'$ in $\{1,..,k\}$. We first show \eqref{convergence_slow_M1}. Since $\lambda_n=\lambda n^\gamma$, multiplying \eqref{m1_sec:slow} by $n^{\alpha-\gamma}$ yields for $j=1,...,k$
\begin{equation}\label{m1_sec:slow_proof}
n^{\alpha-\gamma} M_j^{(n)}(t)= \lambda e^{\lambda t (P-I)}\int_0^t n^\alpha \PP \left(L_j^{(n)}>v\right) e^{-\lambda v (P-I)} \Delta_j P_n e^{\lambda v (P-I)} dv .
\end{equation}
By definition of $L_j^{(n)}$ and the fat tail property of $L_j$:
$$n^\alpha \PP \left(L_j^{(n)}>v\right)=n^\alpha \PP \left(L_j>nv\right)\sim n^\alpha \frac{1}{(nv)^\alpha}=\frac{1}{v^\alpha},\quad v\in (0,t),\quad n\to \infty.$$
Now, since $\lim_{n\to\infty}P_n=P$ and
\begin{multline}\label{m1_sec_argument}
\sup_{n\in \nbN}n^\alpha \PP \left(L_j^{(n)}>v\right)=\sup_{n\in \nbN}n^\alpha \PP \left(L_j>nv\right)=\frac{\sup_{n\in \nbN}(nv)^\alpha \PP \left(L_j>nv\right)}{v^\alpha}\\
\le \frac{\sup_{u\ge 0}u^\alpha \PP \left(L_j>u\right)}{v^\alpha},\quad v\in(0,1),
\end{multline}
the dominated convergence enables us to let $n\to\infty$ in \eqref{m1_sec:slow_proof} to get \eqref{convergence_slow_M1}. We now turn to \eqref{convergence_slow_M20}. Multiplying \eqref{m2_jj_sec:slow} by $n^{\alpha-\gamma}$ yields
\begin{multline}\label{m2_jj_sec:slow_proof}
n^{\alpha-\gamma} M_{jj}^{(n)}(t) = \lambda e^{\lambda t (P-I)}\int_0^t n^\alpha \PP\left(L_j^{(n)}>v\right) e^{-\lambda v (P-I)}
\Delta_j^2 P_n e^{\lambda v (P-I)}dv\\
 +  2 \lambda e^{\lambda t (P-I)}\int_0^t n^\alpha \PP\left(L_j^{(n)}>v\right) \Delta_j P_n M_{j}^{(n)}(v)dv.
\end{multline}
Since \eqref{convergence_slow_M1} in particular implies that $\lim_{n\to \infty}M_{j}^{(n)}(v)=0$ for all $v\in(0,1)$, and thanks to the upper bound \eqref{m1_sec_argument}, a dominated convergence argument entails that the second integral on the right hand side of \eqref{m2_jj_sec:slow_proof} tends to $0$ as $n\to\infty$. We also conclude by a dominated convergence argument that the first integral on the right hand side of \eqref{m2_jj_sec:slow_proof} tends to the right hand side of \eqref{convergence_slow_M20}, and we are done. As to \eqref{m2_jj'_sec:slow}, we have for $j\neq j'$
\begin{multline}\label{m2_jj'_sec:slow_proof}
n^{\alpha-\gamma}M_{jj'}^{(n)}(t) = \lambda e^{\lambda t (P-I)}\int_0^t n^\alpha \PP\left(L_j^{(n)}>v\right)\PP\left(L_{j'}^{(n)}>v\right) e^{-\lambda v (P-I)}
\Delta_j\Delta_{j'} P_n e^{\lambda v (P-I)} dv\\
 + \lambda e^{\lambda t (P-I)}\int_0^t \left\{n^\alpha \PP\left(L_j^{(n)}>v\right) \Delta_j P_n M_{j'}^{(n)}(v)
 +  n^\alpha\PP\left(L_{j'}^{(n)}>v\right) \Delta_{j'} P_n M_{j}^{(n)}(v)\right\} dv.
\end{multline}
Similarly to the second integral on the right hand side of \eqref{m2_jj_sec:slow_proof}, we show that the second integral on the right hand side of \eqref{m2_jj'_sec:slow_proof} converges to $0$ as $n\to\infty$. As to the first integral, the fact that $\PP\left(L_{j'}^{(n)}>v\right)=\PP\left(L_{j'}>nv\right)\longrightarrow 0$ as $n\to\infty$, combined with the upper bound \eqref{m1_sec_argument}, yields by the dominated convergence theorem that it tends to $0$ as $n\to\infty$, achieving the proof of \eqref{m2_jj'_sec:slow} and of the theorem.

\section{A remark on the computation of the limiting joint Laplace transform when $\alpha\in \nbQ$}\label{sec:remark_compute}
We identified in Theorem \ref{theo_regimes} the different limiting regimes when $\gamma$ is larger or equal to $\alpha$ by obtaining the corresponding limiting joint Laplace transform $\chi(s,t)$ in each case. Even though the distributional limits \eqref{convergence_fast} and \eqref{convergence_equilibrium} involve simple processes $\{{\cal X}^\alpha (t),\ t\in[0,1] \}$ and $\{\nu_j^\alpha(t),\ t\ge 0 \}$, $j=1,...,k$, it turns out that the Laplace transforms $\chi(s,t)$, which are solutions to the differential equations \eqref{Poisson_ODE_fast} and \eqref{Poisson_ODE_equilibrium}, are in general not explicit in the fast or equilibrium arriving cases. We suggest to show that things are much simpler when $\alpha\in(0,1)$ is rational, say of the form $$\alpha=1-p/q$$ for some $p$ and $q\in\nbN^*$, with $p<q$. The idea here is quite simple and standard, and consists in expanding a transformation of the solution $t\in[0,1]\mapsto \chi(s,t)\in \nbR^{{\cal S}\times {\cal S}}$ into a power series with matrix coefficients, as explained in \cite[Section 1.1]{Balser00}. Let us focus on the fast arrival case in Section \ref{sec:fast}, although the method is of course applicable to the equilibrium case, and let us put $\check\chi(s,t):= \chi(s,t^p)$, $t\in[0,1]$. In that case, we deduce from \eqref{Poisson_ODE_fast} that $t\in[0,1]\mapsto \check\chi(s,t)$ verifies the matrix differential equation
$$\left\{
\begin{array}{rcl}
\partial_t \check\chi(s,t) &=& \left[(p+q) t^{q}\lambda (P-I) + pt^{p-1}\beta \lambda \sum_{\ell =1}^k s_\ell\Delta_\ell \right]\check\chi(s,t) ,\quad t\in [0,1],\\
&=& [Q_1 t^q + Q_2(s) t^{p-1}]\check\chi(s,t),\\
\chi(s,0)&=& I.
\end{array}
\right.
$$
where $Q_1:=(p+q)\lambda (P-I)$ and $Q_2(s):=p\beta \lambda \sum_{\ell =1}^k s_\ell\Delta_\ell$, $s=(s_1,...,s_k)$. It is quite simple to check that $\check\chi(s,t)$ can then be expanded as
\begin{equation}\label{expansion_chi}
\check\chi(s,t)=\sum_{j=0}^\infty U_j(s) t^j,\quad t\in [0,1],
\end{equation}
where the sequence of matrices $ (U_j(s))_{j\in\nbN}$ is defined from \cite[Relation (1.4)]{Balser00} by $U_0(s)=I$ and
\begin{equation}\label{rel_U_j}
U_j(s)=
\left\{
\begin{array}{c l}
0, & 1\le j <p,\\
Q_2(s)U_{j-p}(s)/j, &  p\le j<q+1,\\
\left[ Q_2(s)U_{j-p}(s)+Q_1 U_{j-q-1}(s) \right]/j, & j\ge q+1 ,
\end{array}
\right.
\end{equation}
and that \eqref{expansion_chi} converges for all $t$, as proved in \cite[Lemma 1 p.2]{Balser00}. The final solution is then expressed in that case as
$$
\chi(s,t)=\check\chi(s,t^{1/p})=\sum_{j=0}^\infty U_j(s) t^{j/p},\quad t\in [0,1].
$$
The $U_j(s)$'s, $j\in\nbN$, being simply expressed with the simple linear recurrence \eqref{rel_U_j}, this expansion for $\chi(s,t)$ is then easy to handle as it can be e.g. approximated by truncation.


\bibliographystyle{alpha}



\end{document}